\documentclass[preprint,sort&compress,final,11pt]{elsarticle}

\usepackage{amsmath}
\usepackage{amssymb}
\usepackage{graphicx}
\usepackage{latexsym}
\usepackage{epstopdf}
\usepackage{natbib}
\usepackage{color}
\usepackage{hyperref}

\newtheorem{theorem}{Theorem}[section]

\newtheorem{lemma}[theorem]{Lemma}

\newtheorem{remark}[theorem]{Remark}
\numberwithin{equation}{section}

\def\Proof{\noindent{\bf Proof.}~}
\def\qed{\hfill$\square$\smallskip}

\journal{\empty}
\date{}

\begin{document}

\begin{frontmatter}

\title{Entire solutions to reaction diffusion equations}

\author[au1,au2]{Yang Wang}

\ead[au1]{ywang@mail.bnu.edu.cn}

\author[au1]{Xiong Li\footnote{ Partially supported by the NSFC (11571041)and the Fundamental Research Funds for the Central Universities. Corresponding author.}}

\address[au1]{School of Mathematical Sciences, Beijing Normal University, Beijing 100875, P.R. China.}

\ead[au1]{xli@bnu.edu.cn}

\address[au2]{ School of Mathematical Sciences, Shanxi University, Taiyuan, Shanxi 030006, P.R. China.}

\begin{abstract}
In this paper, we first use the super-sub solution method to prove the local exponential asymptotic stability of some entire solutions to reaction diffusion equations, including the bistable and monostable cases. In the bistable case, we not only obtain the similar asymptotic  stability result given by Yagisita in 2003, but also simplify his proof. For the monostable case, it is the first time to discuss the local asymptotic stability of entire solutions. Next, we will discuss the asymptotic behavior of entire solutions of bistable equations as $t\rightarrow+\infty$, since the other side was completely known. Here, our results are obtained by use of the asymptotic stability of constant solutions and pairs of diverging traveling front solutions of these equations, instead of constructing the corresponding super-sub solutions as usual.
\end{abstract}

\begin{keyword}
Entire solutions, Traveling front solutions, Reaction diffusion equations
\end{keyword}

\end{frontmatter}

\section{Introduction}
In this paper we focus on the following reaction diffusion equation
\begin{equation}\label{eq:bi}
\partial_tu=\partial_{xx}u+f(u),\ \ \ \ \ \ x\in\mathbb{R},
\end{equation}
where the reaction term $f$ satisfies\\
(A)\ $f\in C^2(\mathbb{R})$, $f(0)=f(\alpha)=f(1)=0$ and $\alpha$ is the unique zero point of $f$ in the interval $(0,1)$, $f'(0)$, $f'(1)<0$; or\\
(A$'$)\ $f\in C^2(\mathbb{R})$, $f(0)=f(1)=0$ and $f(u)>0$ for $u\in(0,1)$, $f'(0)>0$, $f'(1)<0$.

Under the assumption (A), \eqref{eq:bi} is a bistable equation and the background can be found in \cite{aw78}, \cite{c71}, \cite{nya65} and the references therein. This model can illustrate that a nerve has been treated with certain toxins as stated in \cite{c71}. Also it can be used to describe a bistable active transmission line introduced in \cite{nya65}. For more general reaction terms, Aronson and Weinberger in \cite{aw75} used it to describe the heterozygote inferiority case and also pointed out that some flame propagation problems in chemical reactor theory can be demonstrated by equations of the form \eqref{eq:bi}. While under the assumption (A$'$), also in \cite{aw75}, they used it to describe the heterozygote intermediate
case and at this time, \eqref{eq:bi} becomes the famous KPP-Fisher equation, namely the monostable equation, which had been studied in \cite{f37} and \cite{kpp37}.

In recent years, the existence, uniqueness, stability and other properties of traveling wave solutions of \eqref{eq:bi} have been investigated extensively, for example, see \cite{aw75}, \cite{aw78}, \cite{f79l}, \cite{f79}, \cite{fm77}, \cite{fm81}, \cite{f37}, \cite{kpp37} and the references therein. More precisely, for the bistable case, the existence and uniqueness of traveling wave solutions one can refer to \cite{fm77}, while for the monostable case, one can refer to \cite{aw78}, \cite{kpp37} and the references therein. A function $\phi(\xi)$, $\xi=x+ct$, is called a traveling wave solution of \eqref{eq:bi} connecting $0$ and $1$ with the wave speed $c$, if it satisfies
\begin{equation}\label{eq:obi}
\phi''(\xi)-c\phi'(\xi)+f(\xi)=0,\ \ \lim\limits_{\xi\rightarrow-\infty}\phi(\xi)=0,\ \ \lim\limits_{\xi\rightarrow+\infty}\phi(\xi)=1,
\end{equation}
which is actually monotone increasing proved in \cite{fm77} and \cite{om99d}. The reflect $\phi(-x-ct)$ also admits the monotone decreasing traveling wave solution with the opposite wave speed $-c$. In fact, $\phi(-x-ct)$ satisfies $$\lim\limits_{\tilde{\xi}\rightarrow-\infty}\phi(\tilde{\xi})=1, \ \ \ \ \lim\limits_{\tilde{\xi}\rightarrow+\infty}\phi(\tilde{\xi})=0$$ with $\tilde{\xi}=-x-ct$. Thus, if there is a solution to \eqref{eq:obi}, then a traveling wave solution with an opposite speed exists simultaneously. Moreover, the monotone traveling wave solution is also known as the traveling front solution.

However, it is not enough to understand the dynamical structure of \eqref{eq:bi} by only considering traveling wave solutions. Recently, the existence of entire solutions, which are classical solutions and defined for all $(x,t)\in\mathbb{R}\times\mathbb{R}$, is widely discussed. In \cite{hn99}, under the assumption (A$'$) and $f'(u)\leqslant f'(0)\ \ (u\in[0,1])$, Hamel and Nadirashvili proved the existence of entire solutions by the comparison theorem and super-sub solution method, which consists of traveling front solutions and solutions to the diffusion-free equation. Meanwhile, they also pointed out that the solutions to \eqref{eq:bi} depending only on $t$ and traveling wave solutions are typical examples of entire solutions and showed various entire solutions of \eqref{eq:bi} in their subsequent paper \cite{hn01}. While under the assumption (A) and $\int^1_0f(u)du>0$, which implies that wave speeds of any traveling front solutions of \eqref{eq:bi} must be positive, Yagisita in \cite{y03} revealed that the annihilation process is approximated by a backward global solution of \eqref{eq:bi}, which is actually an entire solution. For Allen-Cahn equation
$$\partial_tu=\partial_{xx}u+u(1-u)(u-a)$$
with $a\in(0,1)$, which is a special example of \eqref{eq:bi}, Fukao, Morita and Ninomiya in \cite{fmn04} proposed a simple proof for the existence of entire solutions, which were already found in \cite{y03} by using the super-sub solution method and the exact traveling front solutions. Moreover, Guo and Morita in \cite{gm05} extended the conclusions in \cite{hn99} and \cite{y03} to more general case. Specially, under the assumption (A) and $\int^1_0f(u)du<0$, which implies that wave speeds $c$ of any traveling front solutions of \eqref{eq:bi} must be negative, Chen and Guo in \cite{cg05} used the quite different method to construct the super-sub solutions to obtain the existence and uniqueness of entire solutions of \eqref{eq:bi}, which are different from those in \cite{fmn04}, \cite{gm05} and \cite{y03}. From the dynamical view the study of entire solutions is essential for a full understanding of the transient dynamics and the structures of the global attractor as mentioned in \cite{mn06}. Other papers about the existence of entire solutions, one can refer to \cite{mt09}, \cite{wl15} and \cite{wlr09}.

In this paper we will firstly investigate the asymptotic behavior of entire solutions of bistable reaction diffusion equations as $t\rightarrow+\infty$, since the authors in \cite{cg05} and \cite{gm05} had obtained the exact asymptotic behavior as $t\rightarrow-\infty$. We conjecture that the long time behavior of entire solutions of \eqref{eq:bi} found in \cite{cg05} and \cite{gm05} may be controlled by some asymptotic stable states of \eqref{eq:bi} defined in \cite{f79l}. Luckily, Fife in \cite{f79l} pointed out several kinds of asymptotic stable states, including constant solutions $u\equiv0$, $u\equiv1$, traveling wave solutions, diverging pairs of traveling wave solutions. With these results, we can obtain the long time behavior of these entire solutions. Secondly, Yagisita in \cite{y03} proved that the entire solution of \eqref{eq:bi} is local exponential asymptotic stable by the asymptotic stability of the constructed invariant manifold. Also the authors in \cite{wlr09} obtained the local Lyapunov stability of entire solutions found in \cite{cg05} and \cite{gm05} by the super-sub solution method. By means of establishing the different super-sub solution of \eqref{eq:bi} from that in \cite{wlr09}, we further obtain the local exponential asymptotic stability of the entire solutions of bistable reaction diffusion equations found in \cite{cg05} and \cite{gm05}. Finally, as far as we know, there is no any results about the stability of entire solutions of \eqref{eq:bi} under the assumption (A$'$). Here, we will prove the local exponential asymptotic stability of the entire solutions of the Fisher-KPP equation found in \cite{gm05}.

The rest is organized as follows. For the reader's convenience, in Section 2 we will prove an interior Schauder estimate that had been stated in \cite{fm77} without the proof. In Section 3, we show some known results about the asymptotic stability of constant solutions, traveling wave solutions and diverging pairs of traveling wave solutions of bistable reaction diffusion equations. Then, in Section 4, under the assumption (A), we investigate the long time behavior of entire solutions of \eqref{eq:bi} and prove the local exponential asymptotic stability of entire solutions by the super-sub solution method. Finally, in Section 5, we will prove the local  exponential asymptotic stability of entire solutions of \eqref{eq:bi} under the assumption (A$'$).

\section{An interior Schauder estimate}

In the sequel, we will investigate the long time behavior of entire solutions of \eqref{eq:bi} in the Banach space
$$C_{unif}:=\{u\in C(\mathbb{R}):u\ \rm{is\ bounded\ and\ uniformly\ continuous\ in\ \mathbb{R}}\}$$
with the norm $\|u\|:=\sup\limits_{x\in\mathbb{R}}|u(x)|$. For this purpose, we introduce the initial condition
\begin{equation}\label{eq:in}
u(x,0)=u_0(x),\ \ \ \ \ \ \ \ \ \ x\in\mathbb{R},
\end{equation}
where $u_0(x)\in C_{unif}(\mathbb{R})$, and let the function $u(x,t;u_0)$ be the solution of equation \eqref{eq:bi} with the initial condition \eqref{eq:in}. Now we introduce an interior Schauder estimate that had been stated in \cite{fm77} without the proof. Although the authors in \cite{wlr09} had proven such estimate for more complicated equations, for the reader's convenience, we use the method in \cite{wlr09} and give the proof.
\begin{lemma}\label{thm:sest}
Suppose that $u(x,t;u_0)$ is a bounded solution to \eqref{eq:bi} with \eqref{eq:in} for $(x,t)\in\mathbb{R}\times[0,+\infty)$, and $\|u(\cdot,t;u_0)\|\leqslant L_0$ for some constant $L_0>0$ and all $t\geq 0$. Assume that $f\in C^2(\mathbb{R})$, and  there exists a positive constant $L_1$ such that $\|f\|$, $\|f'\|$ and $\|f''\|\leqslant L_1$ on $[-L_0,L_0]$. Then there is a positive constant $L$ such that $\|\partial_tu(\cdot,t;u_0)\|$, $\|\partial_xu(\cdot,t;u_0)\|$, $\|\partial_{xx}u(\cdot,t;u_0)\|\leqslant L$ for all $t\in[1,+\infty)$, where $L$ depends on $L_0$ and $L_1$ only.
\end{lemma}
\Proof Fix $r>1$, for $s>0$ and $t\in[s+1,s+r]$, we have
$$u(x,t;u_0)=\int^{+\infty}_{-\infty}\frac{u(y,s;u_0)}
{2\sqrt{\pi(t-s)}}
e^{-\frac{(x-y)^2}{4(t-s)}}dy+\int^{t}_{s}\int^{+\infty}_{-\infty}
\frac{f(u(y,\tau;u_0))}{2\sqrt{\pi(t-\tau)}}
e^{-\frac{(x-y)^2}{4(t-\tau)}}dyd\tau.$$
By the dominated convergence theorem, we get
\begin{align}
\partial_xu(x,t;u_0)=&-\int^{+\infty}_{-\infty}\frac{(x-y)u(y,s;u_0)}
{4\sqrt{\pi}(t-s)^{\frac{3}{2}}}
e^{-\frac{(x-y)^2}{4(t-s)}}dy
\nonumber\\[0.2cm]
&-\int^{t}_{s}\int^{+\infty}_{-\infty}
\frac{(x-y)f(u(y,\tau;u_0))}{4\sqrt{\pi}(t-\tau)^{\frac{3}{2}}}
e^{-\frac{(x-y)^2}{4(t-\tau)}}dyd\tau,\label{eq:est1}
\end{align}
which implies that
\begin{align*}
|\partial_xu(x,t;u_0)|&\leqslant\frac{L_0}{2\sqrt{\pi(t-s)}}
\int^{+\infty}_{-\infty}\frac{|y|}{2(t-s)}e^{-\frac{y^2}
{4(t-s)}}dy\\[0.2cm]
&\ \ \ \ +\frac{L_1}{2\sqrt{\pi}}\int^{t}_{s}\int^{+\infty}_{-\infty}
\frac{|y|}{2(t-\tau)^{\frac{3}{2}}}e^{-\frac{y^2}{4(t-\tau)}}
dyd\tau\\[0.2cm]
&\leqslant\frac{L_0}{\sqrt{\pi}}+\frac{L_1}{\sqrt{\pi}}\int^{t}_{s}
\frac{1}{\sqrt{t-\tau}}d\tau=\frac{L_0}{\sqrt{\pi}}+\frac{2L_1}{\sqrt{\pi}}
(t-s)^{\frac{1}{2}}\\[0.2cm]
&\leqslant\frac{L_0}{\sqrt{\pi}}+\frac{2L_1}{\sqrt{\pi}}
r^{\frac{1}{2}}:=L_2.
\end{align*}
Due to the arbitrariness of $s$, then $\|\partial_xu(\cdot,t;u_0)\|\leqslant L_2$ for all $t\geq 1$.

From \eqref{eq:est1}, one can get
\begin{align*}
\partial_xu(x,t;u_0)=&-\int^{+\infty}_{-\infty}\frac{(x-y)u(y,s;u_0)}
{4\sqrt{\pi}(t-s)^{\frac{3}{2}}}
e^{-\frac{(x-y)^2}{4(t-s)}}dy\\[0.2cm]
&+\int^{t}_{s}\int^{+\infty}_{-\infty}
\frac{f'(u(y,\tau;u_0))\partial_yu(y,s;u_0)}
{2\sqrt{\pi(t-\tau)}}
e^{-\frac{(x-y)^2}{4(t-\tau)}}dyd\tau.
\end{align*}
Thus we have
\begin{align*}
\partial_{xx}u(x,t;u_0)=&-\int^{+\infty}_{-\infty}\frac{u(y,s;u_0)}
{4\sqrt{\pi}(t-s)^{\frac{3}{2}}}
e^{-\frac{(x-y)^2}{4(t-s)}}dy\\[0.2cm]
&+\int^{+\infty}_{-\infty}\frac{(x-y)^2u(y,s;u_0)}
{8\sqrt{\pi}(t-s)^{\frac{5}{2}}}
e^{-\frac{(x-y)^2}{4(t-s)}}dy\\[0.2cm]
&-\int^{t}_{s}\int^{+\infty}_{-\infty}
\frac{(x-y) f'(u(y,\tau;u_0))}
{4\sqrt{\pi}(t-\tau)^{\frac{3}{2}}}
\partial_yu(y,s;u_0)e^{-\frac{(x-y)^2}
{4(t-\tau)}}dyd\tau,
\end{align*}
which also implies that
\begin{align}
|\partial_{xx}u(x,t;u_0)|\leqslant&~\frac{L_0}{2(t-s)}+
L_0\int^{+\infty}_{-\infty}\frac{(x-y)^2}
{8\sqrt{\pi}(t-s)^{\frac{5}{2}}}
e^{-\frac{(x-y)^2}{4(t-s)}}dy\nonumber\\[0.2cm]
&+L_1L_2\int^{t}_{s}\int^{+\infty}_{-\infty}
\frac{|y|}{4\sqrt{\pi}(t-\tau)^{\frac{3}{2}}}
e^{-\frac{y^2}{4(t-\tau)}}dyd\tau\nonumber\\[0.2cm]
\leqslant&~\frac{L_0}{t-s}+2L_1L_2\frac{r^{\frac{1}{2}}}
{\sqrt{\pi}}
\leqslant~L_0+2L_1L_2\frac{r^{\frac{1}{2}}}
{\sqrt{\pi}}:=L_3.\label{eq:est2}
\end{align}
Consequently, $\|\partial_{xx}u(\cdot,t;u_0)\|\leqslant L_3$ for $t\geq 1$.

It follows from \eqref{eq:bi} and \eqref{eq:est2} that $$|\partial_tu(x,t;u_0)|\leqslant|\partial_{xx}u(x,t;u_0)|+
|f(u)|\leqslant L_3+L_1:=L_4.$$
By setting $L=\max\{L_2,L_3,L_4\}$, we finish the proof.
\qed

\section{Some known results about the bistable equation}
In this section and the next section, we consider the bistable case. Now we list the known results about the stability  of constant solutions, traveling wave solutions and diverging pairs of traveling wave solutions (constructed by traveling wave solutions and their reflects), which will be used later.

We firstly state the asymptotic stability of constant solutions $u\equiv0$ and $u\equiv1$, which had been proved in \cite{f79}.
\begin{lemma}\label{lem:stacon}
Suppose that $0\leqslant u_0\leqslant1$ is a continuous function.
\begin{enumerate}
\item If $\int^{1}_{0}f(u)du\geqslant0$, ~$\liminf\limits_{x\rightarrow
      \pm\infty}u_0(x)>\alpha$, then $\lim\limits_{t\rightarrow+\infty}\|u(x,t;u_0)-1\|=0$;
\item if $\inf\limits_{x\in\mathbb{R}}u_0(x)>\alpha$, then $\lim\limits_{t\rightarrow+\infty}\|u(x,t;u_0)-1\|=0$;
\item if $\int^{1}_{0}f(u)du\leqslant0$, ~$\limsup\limits_{x\rightarrow
    \pm\infty}u_0(x)<\alpha$, then $\lim\limits_{t\rightarrow+\infty}\|u(x,t;u_0)\|=0$;
\item if $\sup\limits_{x\in\mathbb{R}}u_0(x)<\alpha$, then $\lim\limits_{t\rightarrow+\infty}\|u(x,t;u_0)\|=0$.
\end{enumerate}
\end{lemma}

Secondly, the global exponential asymptotic stability of traveling wave solutions of \eqref{eq:bi} had been proved in \cite{fm77}.
\begin{lemma}\label{lem:statra}
Suppose that $\phi$ is the solution to \eqref{eq:obi} and $0\leqslant u_0\leqslant1$ is a continuous function.
\begin{enumerate}
  \item If $\limsup\limits_{x\rightarrow-\infty}u_0(x)<\alpha$, $\liminf\limits_{x\rightarrow+\infty}u_0(x)>\alpha$, then there are some constants $x_0$, $M_1>0$ and $\omega_1>0$\ such that
      $$|u(x,t;u_0)-\phi(x+ct-x_0)|<M_1 e^{-\omega_1t},\qquad t\geq 0,\ \  x\in\mathbb{R};$$
  \item if $\limsup\limits_{x\rightarrow+\infty}u_0(x)<\alpha$, $\liminf\limits_{x\rightarrow-\infty}u_0(x)>\alpha$, then there are some constants $x_1$, $M_2>0$ and $\omega_2>0$ such that
      $$|u(x,t;u_0)-\phi(-x+ct-x_1)|<M_2 e^{-\omega_2t},\qquad t\geq 0,\ \ x\in\mathbb{R}.$$
\end{enumerate}
\end{lemma}

In addition, the local asymptotic stability of traveling wave solutions of \eqref{eq:bi} had been proved in \cite{om99d} and \cite{om99p} by the  different ways compared with the method in \cite{fm77}. The stability of traveling wave solutions in $\mathbb{R}^n$ can be found in \cite{lx92}, \cite{r04} and \cite{x92}. For more general reaction terms, the authors in the papers \cite{gr07} and \cite{r08} proved the stability of traveling wave solutions.

Thirdly, as stated in \cite{f79l}, there are four kinds of bounded stationary solutions of equation \eqref{eq:bi}, which are solutions of
\begin{equation}\label{eq:ell}
\frac{d^2 u}{dx^2}+f(u)=0,
\end{equation}
that is, the solution with a single minimum point, the solution with a single maximum point, the periodic solution and the monotone solution. In fact, the solution with a single maximum point or a single minimum point is the homoclinic orbit. More importantly, the author also proved that the solution of \eqref{eq:ell} with a maximum or minimum at a finite value of $x$ are unstable, which means that only the monotone solution, namely the heteroclinic orbit, is stable. From the views of parabolic equations, the monotone solution of \eqref{eq:ell} is the solution to \eqref{eq:obi} with $c=0$ and the stability had been proved in \cite{fm77}. Specially, when $f(u)=u(1-u)(u-a)$ with $a\neq\frac{1}{2}$,  the author in \cite{k96} pointed out that $\phi(\xi)$ is the solution of \eqref{eq:ell} with $\phi(+\infty)=\phi(-\infty)$ if and only if $c=0$. However, in the classical two species Lotka-Volterra competition model, there exists a similar homoclinic orbit with $c\neq0$.

Finally, when $\int^{1}_{0}f(u)du\neq0$, the authors in \cite{fm77} discussed the asymptotic stability of  diverging pairs of traveling wave solutions, where the initial function $u_0(x)$ is assumed to lie above the line $u=\alpha$ for $x$ in the large symmetrical interval about the origin and $\limsup\limits_{x\rightarrow\pm\infty}u_0(x)<\alpha$, or lie below the line $u=\alpha$ for $x$ in the large symmetrical interval about the origin and $\liminf\limits_{x\rightarrow\pm\infty}u_0(x)>\alpha$, that is, there are two pairs of positive constants $\beta_1, \beta_2$, $\overline{L}_1, \overline{L}_2$ such that $u_0(x)>\alpha+\beta_1$ for $|x|<\overline{L}_1$, or $u_0(x)<\alpha-\beta_2$ for $|x|<\overline{L}_2$, where the large constants $\overline{L}_1, \overline{L}_2$ depend on $\beta_1, \beta_2$ and $f$. However, the explicit expressions  of $\overline{L}_1, \overline{L}_2$  are not given in \cite{fm77}. Now, according to the proof of Lemma 6.1 in \cite{fm77}, we can give a low bound of $\overline{L}_1, \overline{L}_2$, respectively.

First we give two estimates. It is easy to see that the characteristic equations of the linearized equations of \eqref{eq:obi} at $u=0$ and $u=1$ are
$$\lambda^2-c\lambda+f'(0)=0,\ \ \ \ \ \ \  \ \mu^2-c\mu+f'(1)=0.$$
The corresponding eigenvalues are
\begin{align*}
\lambda_1=\frac{c+\sqrt{c^2-4f'(0)}}{2},\ \ \ \ \ \ \lambda_2=\frac{c-\sqrt{c^2-4f'(0)}}{2},\\[0.2cm]
\mu_1=\frac{c+\sqrt{c^2-4f'(1)}}{2},\ \ \ \ \ \  \mu_2=\frac{c-\sqrt{c^2-4f'(1)}}{2}.
\end{align*}
Therefore, there are some positive constants $M_3$, $\widetilde{M}_3$, $M_4$ and $\widetilde{M}_4$  such that
\begin{equation}\label{eq:esti}
\begin{array}{ll}
\widetilde{M}_3e^{\mu_2\xi}\leqslant1-\phi(\xi)\leqslant M_3e^{\mu_2\xi},\qquad&\xi\geqslant0,
\\[0.2cm]
\widetilde{M}_4e^{\lambda_1\xi}\leqslant\phi(\xi)\leqslant M_4e^{\lambda_1\xi},&\xi\leqslant0,
\end{array}
\end{equation}
which are also given in \cite{gm05}.

Since $f'(\alpha)>0$, then
\begin{equation}\label{eq:con}
w:=\max\limits_{u\in[0,1]}f'(u)>0.
\end{equation}
Also since $\lim\limits_{u\rightarrow1^-}\frac{f(u)}{1-u}=-f'(1)$ and $\lim\limits_{u\rightarrow0^+}\frac{f(u)}{u}=f'(0)$, then the functions $\frac{f(u)}{1-u}$ and $\frac{f(u)}{u}$ are continuous in the interval $[0,1]$ if we define the function $\frac{f(u)}{1-u}$ as $-f'(1)$ at $u=1$,  $\frac{f(u)}{u}$ as $f'(0)$ at $u=0$. Thus there is a positive constant $b$ such that $$
|f(u)|\leqslant b(1-u),\qquad|f(u)|\leqslant bu,\qquad u\in[0,1].
$$

Give any $\beta_1>0$, choose $q_0$, $q_1$ as $0<1-q_1<1-q_0<\alpha+\beta_1$, let $\tilde{\mu}_1>0$, $\tilde{\beta}$ and $\widetilde{M}$ be corresponding to $\mu_1$, $\beta$ and $M$ in the proof of Lemma 6.1 in \cite{fm77} respectively, choose $\tilde{\mu}_2, \overline{M}, \varphi_0$ as
\begin{align*}
0&<\tilde{\mu}_2<\min\{-\mu_2c,\ \tilde{\mu}_1\},\ \ \overline{M}=\frac{w+b}{c\beta\mu_2}M_3-\frac{w+\tilde{\mu}_2}
{\beta\tilde{\mu}_2}q_0<0, \\[0.2cm]
\varphi_0&<\min\left\{\overline{M},\ \overline{M}-\frac{1}{\mu_2}
\ln\frac{q_1-q_0}{M_3},\ \overline{M}-\frac{1}{\mu_2}
\ln\frac{(\tilde{\mu}_1-\tilde{\mu}_2)q_0}{bM_3}\right\},
\end{align*}
then we can obtain a low bound of $\overline{L}_1$ as
$$
\overline{L}_1\geqslant \widetilde{M}\geqslant\max\left\{-\varphi_0,\ -\frac{1}{\lambda_1}\ln
\frac{1-\alpha-\beta_1}{M_4}-\varphi_0\right\}.
$$
Similarly, give any $0<\beta_2<\alpha$, choose $\tilde{q}_0$, $\tilde{q}_1$ as $0<\alpha-\beta_2<\tilde{q}_0<\tilde{q}_1<\alpha$, let $\tilde{\mu}'_1>0$, $\widetilde{M}'$ and $\tilde{\beta}'$  be corresponding to $\mu_1'$, $M'$ and $\beta'$ in the proof of Lemma 6.1 in \cite{fm77} respectively, choose $\tilde{\mu}'_2, \overline{M}', \tilde{\varphi}_0$ as
\begin{align*}
 0&<\tilde{\mu}'_2<\min\{-\lambda_1c,\ \tilde{\mu}'_1\},\ \ \overline{M}'=\frac{w+b}{c\lambda_1\beta'}M_4-\frac{w+\tilde{\mu}'_2}
{\beta\tilde{\mu}'_2}\tilde{q}_0<0,\\[0.2cm]
\tilde{\varphi}_0&<\min\left\{\overline{M}',\ \overline{M}'+
\frac{1}{\lambda_1}\ln\frac{\tilde{q}_1-\tilde{q}_0}{M_4},
\ \overline{M}'+\frac{1}{\lambda_1}
\ln\frac{(\tilde{\mu}'_1-\tilde{\mu}'_2)\tilde{q}_0}{bM_4}\right\},
\end{align*}
then we can obtain a low bound of $\overline{L}_2$ as
$$
\overline{L}_2\geqslant \widetilde{M}'\geqslant\max\left\{-\tilde{\varphi}_0,\ \frac{1}{\mu_2}\ln
\frac{\alpha-\beta_2}{M_3}-\tilde{\varphi}_0\right\}.
$$

Therefore we can refine Theorem 3.2 in \cite{fm77} and obtain the following lemma.
\begin{lemma}\label{lem:twd}
Suppose that $\phi$ is the solution to \eqref{eq:obi}, and $0\leqslant u_0\leqslant1$ is a continuous function.
\begin{enumerate}
\item If $\int^{1}_{0}f(s)ds>0$, $\limsup\limits_{x\rightarrow\pm\infty}u_0(x)<\alpha$ and $u_0(x)>\alpha+\beta_1$ for $|x|<\overline{L}_1$,
 then there are some constants $x_2$, $x_3$ and positive constants $M_5$, $\omega_3$ such that
\begin{align*}
&|u(x,t;u_0)-\phi(x+ct-x_2)|<M_5 e^{-\omega_3t},\ \ \qquad x<0,\ \ t\geq 0;\\[0.2cm]
&|u(x,t;u_0)-\phi(-x+ct-x_3)|<M_5 e^{-\omega_3t},\qquad x>0,\ \ t\geq 0.
\end{align*}
\item If $\int^{1}_{0}f(s)ds<0$, $\liminf\limits_{x\rightarrow\pm\infty}u_0(x)>\alpha$ and $u_0(x)<\alpha-\beta_2$ for $|x|<\overline{L}_2$,
then there are some constants $x_4$, $x_5$ and positive constants $M_6$, $\omega_4$ such that
\begin{align*}
&|u(x,t;u_0)-\phi(x+ct-x_4)|<M_6 e^{-\omega_4t},\ \ \qquad x>0,\ \ t\geq 0;\\[0.2cm]
&|u(x,t;u_0)-\phi(-x+ct-x_5)|<M_6 e^{-\omega_4t},\qquad x<0,\ \ t\geq 0.
\end{align*}
\end{enumerate}
\end{lemma}

More interestingly, the author in \cite{f79l} conjectured that when the bounded initial function $u_0$ is away from $\alpha$ for large $|x|$, the asymptotic stable solution may only be one of the four kinds: $u\equiv0$, $u\equiv1$, the traveling wave solution, the diverging pairs of traveling wave solutions, while other solutions are unstable. Although the conjecture is partially solved in the paper \cite{f79}, it has been not completely solved and is still open. Moreover, we also remark that \eqref{eq:bi} admits a traveling wave solution $\hat{\phi}$ connecting $\alpha$ and $1$, and there are some results about $u(x,t;u_0)$ converging to $\hat{\phi}$ when $\alpha\leq u_0\leq 1$, referring to \cite{kpp37}. The authors in \cite{mn06} also presented that if $u_0(x)<\alpha$ for $x$ in some interval of $x$-axis, then wether $u(x,t;u_0)$ converges to $\hat{\phi}$ or not? For example, if $u_0$ satisfies the following condition
\begin{align*}\begin{array}{ll}
\lim\limits_{x\rightarrow+\infty}u_0(x)=1,\
\lim\limits_{x\rightarrow-\infty}u_0(x)=\alpha,\
\mbox{and there exists}~\tilde{x}~\mbox{such that}\ \mbox{when}\\[0.2cm]
x\leqslant\tilde{x},
u_0(x)<\alpha,
\end{array}
\end{align*}
then wether $u(x,t;u_0)$ converges to $\hat{\phi}$ or not? Thus this open problem contains the opposite side of the above conjecture. That is, suppose that the initial condition $u_0(x)$ is bounded and, when $|x|$ is sufficiently large, $u_0(x)$ is not far away from $\alpha$ at least on one side of $x$-axis, then wether $u(x,t;u_0)$ converges to $\hat{\phi}$ or not? In a word, the proof of this open problem will enforce the understanding of the conjecture.

\section{The bistable equation}
In this section, we will discuss the local exponential asymptotic stability of entire solutions of \eqref{eq:bi} found in \cite{cg05} and \cite{gm05} respectively, and their asymptotic behaviors when $t$ converges to $+\infty$.

Before stating the main results, we do some preparations. Firstly, in order to prove the existence of entire solutions, the authors in \cite{gm05} constructed two pairs of different super-sub solutions corresponding to $f'(0)\leqslant f'(1)$ and $f'(0)>f'(1)$, respectively. Moreover, it follows from \cite{fm77} that
$$\rm{if}\hspace{0.3cm}\int^{1}_{0}f(u)du\gtreqqless0,
\hspace{0.5cm}\rm{then}\hspace{0.3cm} c\gtreqqless0.$$
Hence we will discuss the long time behaviors and stabilities of entire solutions under the following four cases:\\[0.2cm]
(C1) $\int^{1}_{0}f(u)du>0$ and $f'(0)>f'(1)$;\\[0.2cm]
(C2) $\int^{1}_{0}f(u)du>0$ and $f'(0)\leqslant f'(1)$;\\[0.2cm]
(C3) $\int^{1}_{0}f(u)du<0$ and $f'(0)>f'(1)$;\\[0.2cm]
(C4) $\int^{1}_{0}f(u)du<0$ and $f'(0)\leqslant f'(1)$.\\[0.2cm]

Secondly, for the sake of the proof of the uniqueness of entire solutions, similar to \cite{cg05}, we introduce the metable dynamics of \eqref{eq:bi}. We call the solution $u(x,t)$ of \eqref{eq:bi} satisfying the condition
$\mathbb{M}^+$, if there exist the constants $d_1>0$ and $T_1\in\mathbb{R}$, the functions $l_1(t)$ and $m_1(t)$ such that for all $t\leqslant T_1$,
$$\left\{\begin{array}{ll}
u(x,t)\leqslant\alpha_1,\ \ \ \forall x\in[\min\{l_1(t)+d_1,m_1(t)-d_1\},\max\{l_1(t)+d_1,
m_1(t)-d_1\}],\\[0.2cm]
u(x,t)\geqslant\alpha_2,\ \ \ \forall x\in(-\infty,l_1(t)]\cup[m_1(t),+\infty),
\end{array}\right.$$
where $\alpha_1$ and $\alpha_2$ are some constants satisfying $f\neq0$ in $(0,\alpha_1]\cup[\alpha_2,1)$. Similarly, the solution $u(x,t)$ of \eqref{eq:bi} is called to satisfy the condition $\mathbb{M}^-$, if there exist the constants $d_2>0$ and $T_2\in\mathbb{R}$, the functions $l_2(t)$ and $m_2(t)$ such that for all $t\leqslant T_2$,
$$\left\{\begin{array}{ll}
u(x,t)\leqslant\alpha_1,\ \forall x\in(-\infty,l_2(t)]\cup[m_2(t),+\infty],\\[0.2cm]
u(x,t)\geqslant\alpha_2,\ \forall x\in[\min\{l_2(t)+d_2,m_2(t)-d_2\},\max\{l_2(t)+d_2,m_2(t)-d_2\}],
\end{array}\right.$$
where $\alpha_1$ and $\alpha_2$ are some constants satisfying $f\neq0$ in $(0,\alpha_1]\cup[\alpha_2,1)$.

We initially discuss entire solutions given in Theorem 1.1 from \cite{gm05}. By the method of the proof of the uniqueness in \cite{cg05} and Lemma \ref{lem:stacon}, we can prove the uniqueness of entire solutions and obtain the following result.
\begin{theorem}\label{thm:entire1}
Suppose that $\int^1_0f(u)du>0$. Let $\phi$ be the solution of \eqref{eq:obi} with the wave speed $c$. Then for any given constants $y_1$ and $y_2$, there is a unique entire solution (up to a translation in $t$ and $x$) $u_1(x,t)$ of \eqref{eq:bi} defined for all $(x,t)\in\mathbb{R}\times\mathbb{R}$ such that $0<u_1(x,t)<1$, $\partial_t u_1(x,t)>0$ and
\begin{equation}\label{eq:asy}
\lim\limits_{t\rightarrow-\infty}\{\sup\limits_{x\geqslant0}
|u_1(x,t)-\phi(x+ct+y_1)|
+\sup\limits_{x\leqslant0}|u_1(x,t)-\phi(-x+ct+y_2)|\}=0,
\end{equation}
$$\lim\limits_{t\rightarrow+\infty}\sup\limits_{x\in\mathbb{R}}|u_1(x,t)-1|=0.$$
\end{theorem}
\Proof We first consider the case (C1). When $t\leqslant0$, it follows from \cite{gm05} that the supersolution and subsolution of \eqref{eq:obi} are $$\overline{u}_1(x,t)=\min\{\phi(x+p_1(t))+\phi(-x+p_1(t)),1\},$$
$$\underline{u}_1(x,t)=\max\{\phi(x+ct+x_6),\phi(-x+ct+x_6)\},$$
where the function $p_1(t)$ $(t\leqslant0)$ is the solution to
\begin{equation}\label{eq:ode1}\left\{\begin{array}{ll}
p'_1(t)=c+M_7e^{\lambda_1p_1(t)},\quad t<0,\\[0.2cm]
p_1(0)\leqslant0,
\end{array}
\right.
\end{equation}
$c$ is the wave speed, $\lambda_1=\frac{c+\sqrt{c^2-4f'(0)}}{2}$, and the expression of $M_7>0$ is too long, which can be found in\cite{gm05}, and
$x_6=p_1(0)-\frac{1}{\lambda_1}\ln(1+\frac{M_7}{c})$.

Substituting $t=0$ into the expression of the subsolution $\underline{u}_1(x,t)$, we have
$$\liminf\limits_{x\rightarrow\pm\infty}u_1(x,0)\geqslant
\liminf\limits_{x\rightarrow\pm\infty}\underline{u}_1(x,0)=
\lim\limits_{x\rightarrow\pm\infty}\max\{\phi(x+x_6),\phi(-x+x_6)\}
=1>\alpha.$$
Also since $\int^{1}_{0}f(s)ds>0$, it follows from Lemma \ref{lem:stacon} that
$$\lim\limits_{t\rightarrow+\infty}\|u_1(x,t)-1\|=0,$$
which is the asymptotic behaviors of the entire solution $u_1(x,t)$ when $t$ converges to $+\infty$. On the other hand, \eqref{eq:asy} has been proved in \cite{mt09}.

Since $\underline{u}_1(x,t+t_1)<\overline{u}_1(x,t)$ for arbitrary $t_1$, the method in \cite{gm05} to prove the uniqueness of entire solutions is no longer valid. Here we use the method in \cite{cg05} to prove the uniqueness, and only need to verify that $u_1(x,t)$ satisfies the condition $\mathbb{M}^+$.

Since $f$ satisfies the assumption (A), there are some constants $\alpha_1<\alpha<\alpha_2$ such that $f\neq0$ in $(0,\alpha_1]\cup[\alpha_2,1)$. Then from the monotonicity of $\phi$, there exists a positive constant $\tilde{x}$ such that $\phi(\tilde{x}+x_6)\geqslant\alpha_2$.
Set $$l_1(t)=ct-\tilde{x}, \ \ \ \ m_1(t)=-ct+\tilde{x}.$$ For any $t\leqslant0$, obviously $m_1(t)\geqslant0\geqslant l_1(t)$. On one hand, since $$u_1(x,t)\geqslant\max\{\phi(x+ct+x_6),\phi(-x+ct+x_6)\},$$
then $u_1(x,t)\geqslant\alpha_2$ for any $x\in(-\infty,l_1(t)]\cup[m_1(t),+\infty)$. On the other hand, according to \eqref{eq:asy}, for any $\varepsilon>0$, there exists a $T_{11}<0$ such that for any $t\leqslant T_{11}$
$$\sup\limits_{x\geqslant0}|u_1(x,t)-\phi(x+ct+y_1)|
\leqslant\varepsilon,\quad\sup\limits_{x\leqslant0}|u_1(x,t)-\phi(-x+ct+y_2)|
\leqslant\varepsilon,$$
which together with the monotonicity of $\phi$ implies that
$$u_1(x,t)\leqslant \phi(m_1(t)-d_1+ct+y_1)+\varepsilon\leqslant\alpha_1$$
for $x\leqslant m_1(t)-d_1=-ct+\tilde{x}-d_1$ with the sufficiently large positive constant $d_1$. Similarly, $u_1(x,t)\leqslant\alpha_1$ for $x\geqslant l_1(t)+d_1$. Since there is a $T_{12}<0$ such that $l_1(t)+d_1\leqslant0\leqslant m_1(t)-d_1$ for any $t\leq T_{12}$, choose  $T_1\leqslant\min\{T_{11},T_{12}\}$, then $u_1(x,t)\leqslant\alpha_1$ for any $t\leqslant T_1$ and $x\in[l_1(t)+d_1,m_1(t)-d_1]$. Therefore, the entire solution $u_1(x,t)$ satisfies the condition $\mathbb{M}^+$, and from the proof in the papers \cite{cg05} and \cite{wlr09}, it is unique except for a space-time translation.

Finally, we will prove the monotonicity of $u_1(x,t)$ with respect to $t$. From the proof of the existence of entire solutions in the papers \cite{cg05} and \cite{gm05} and so on, the authors choose the function $u_n(x,t)$ as the unique classical solution to the following initial problem
\begin{equation}\label{eq:cauchy}
\left\{\begin{array}{ll}
\partial_tu_n=\partial_{xx}u_n+f(u_n),\quad x\in\mathbb{R},\quad t>-n,\\[0.2cm]
u_n(x,-n)=\underline{u}_1(x,-n),\quad x\in\mathbb{R}.
\end{array}
\right.
\end{equation}
Then by Lemma \ref{thm:sest} and the process of diagonalization, there exists a subsequence of $\{u_n\}$ converging in the space $C^{2,1}_{loc}(\mathbb{R}\times(-\infty,0])$. In fact, the limit of this subsequence is the entire solution as we desired. Obviously, $\underline{u}_1(x,t)$ is the subsolution to \eqref{eq:cauchy} and $\partial_t\underline{u}_1(x,t)\geqslant0$ for $(x,t)\in\mathbb{R}\times(-\infty,0]$. It is easy to see that $u_n(x,t)$ satisfies
$$\partial_tu_n|_{t=-n}=\partial_{xx}u_n+f(u_n)
=\partial_{xx}\underline{u}_1+f(\underline{u}_1)\geqslant
\partial_t\underline{u}_1\geqslant0,$$
and since $\partial_tu_n$ satisfies the equation
$\partial_t(\partial_tu_n)-\partial_{xx}(\partial_tu_n)-f'(u_n)\partial_tu_n=0$, then $\partial_tu_n\geqslant0$ for $(x,t)\in\mathbb{R}\times[-n,0]$ by the minimum theorem. From the process of convergence mentioned in the above and the strong minimum theorem, $\partial_tu_1(x,t)>0$ for all $(x,t)\in\mathbb{R}\times\mathbb{R}$.

Now we deal with the case (C2). According to \cite{gm05}, $u_1(x,t)$ is unique and satisfies \eqref{eq:asy}, and the supersolution and subsolution for $t\leqslant0$ are
$$\overline{u}_1(x,t)=\phi(x+p_2(t))+\phi(-x+p_2(t)),\ \ \underline{u}_1(x,t)=\phi(x+p_3(t))+\phi(-x+p_3(t)),$$
where the functions $p_2(t)$ and $p_3(t)$ satisfy
\begin{equation}\label{eq:2ode}
\left\{\begin{array}{ll}
p'_2(t)=c+M_8e^{\lambda_1p_2(t)},\\[0.2cm]
p_2(0)\leqslant0,
\end{array}\right.\ \ \ \ \ \ \ \ \
\left\{\begin{array}{ll}
p'_3(t)=c-M_8e^{\lambda_1p_3(t)},\\[0.2cm]
p_3(0)\leqslant\min\{0,\frac{1}{\lambda_1}\ln(\frac{c}{M_8})\},
\end{array}\right.
\end{equation}
and the expression of $M_8>0$ can be found in \cite{gm05}.

It is easy to see that
$$\liminf\limits_{x\rightarrow\pm\infty}\underline{u}_1(x,0)=
\lim\limits_{x\rightarrow\pm\infty}\{\phi(x+p_3(0))+\phi(-x+p_3(0))\}
=1>\alpha,$$
which together with Lemma \ref{lem:stacon} implies that
$$\lim\limits_{t\rightarrow+\infty}\|u_1(x,t)-1\|=0.$$

Finally we only need to prove the monotonicity of $u_1(x,t)$ with respect to $t$. For this purpose, we first show that $p'_3(t)>0$ for $t\in(-\infty,0]$. Since $p_3(0)\leqslant\min\{0,\frac{1}{\lambda_1}\ln(\frac{c}{M_8})\}$, then $p'_3(0)=c-M_8e^{\lambda p_3(0)}>0$, and if there is a $t_2\in(-\infty,0)$ such that $p'_3(t_2)=0$ and $p'_3(t)>0$ for $t\in(t_2,0)$, then $$\frac{1}{\lambda_1}\ln\frac{c}{M_8}=p_3(t_2)\leqslant p_3(0)<\frac{1}{\lambda_1}\ln\frac{c}{M_8},$$ which is a contradiction. Thus $p'_3(t)>0$ for $t\in(-\infty,0]$, which implies that $$\partial_t\underline{u}_1(x,t)=\phi'(x+p_3(t))p'_3(t)+
\phi'(-x+p_3(t))p'_3(t)>0.$$
Therefore, $\partial_t u_1(x,t)>0$ for all $(x,t)\in\mathbb{R}\times\mathbb{R}$.
\qed

Now we will consider the stability of the unique entire solution $u_1(x,t)$ obtained in Theorem \ref{thm:entire1} by using the method in \cite{wlr09}, and obtain the following result.
\begin{theorem}\label{thm:entire1l}
Suppose that $\int^{1}_{0}f(u)du>0$. Then the unique entire solution $u_1(x,t)$ obtained in Theorem \ref{thm:entire1} is local exponential asymptotic stable.
\end{theorem}
\Proof To begin with, we introduce some notations. Due to $f\in C^2(\mathbb{R})$ and $f'(0)$, $f'(1)<0$, there exists a $\theta>0$ such that $f'(u)<0$ on $[-\theta,2\theta]\cup[1-2\theta,1+\theta]$. Therefore,
\begin{equation}\label{eq:const}
v:=\max\left\{\max\limits_{[-\theta,2\theta]}f'(u),
\max\limits_{[1-2\theta,1+\theta]}f'(u)
\right\}<0.
\end{equation}
Let $\phi(\xi)$ be the solution to \eqref{eq:obi}. When $\phi(x+ct)$ or $\phi(-x+ct)\in[\theta,1-\theta]$, there is a $\beta>0$ such that $$\phi(x+ct)+\phi(-x+ct)\geqslant\beta.$$
In the sequel of this paper, we still use these notations.

We first consider the case (C1). For $t\geqslant0$, we will prove that the following functions
\begin{align*}
\overline{u}_1^{+}(x,t)=\min\{1,u_1(x,t+\gamma(t))+q(t)\},\\[0.2cm] \underline{u}_1^{-}(x,t)=\max\{0,u_1(x,t-\gamma(t))-q(t)\}
\end{align*}
are the supersolution and subsolution to \eqref{eq:bi} with the initial condition $u_0(x)=u_1(x,0)$ respectively, where the function $q(t)$ is the solution to the following initial problem
\begin{equation}\label{eq:ode1}\left\{\begin{array}{ll}
q'(t)-vq(t)=0,\quad t>0,\\[0.2cm]
q(0)=\overline{q}_0,
\end{array}
\right.
\end{equation}
$0\leqslant\overline{q}_0\leqslant\theta$ may be arbitrary, and the function $\gamma(t)$ is to be determined later.

Here, we mainly prove that the function $\overline{u}_1^{+}(x,t)$ is the supersolution, and the rest is similar. When $\overline{u}_1^{+}(x,t)\equiv1$, the conclusion is obvious. Thus, we only consider $\overline{u}_1^{+}(x,t)=u_1(x,t+\gamma(t))+q(t)$. Firstly, when $u_1(x,t)\in[0,\theta]$ or $[1-\theta,1]$, from \eqref{eq:const}, \eqref{eq:ode1} and $\partial_tu_1(x,t)>0$, we can get
\begin{align*}
\partial_t\overline{u}_1^{+}-\partial_{xx}\overline{u}_1^{+}
-f(\overline{u}_1^{+})&=\gamma'(t)
\partial_tu_1+\partial_tu_1-\partial_{xx}u_1+q'(t)-f(u_1+q(t))\\[0.2cm]
&\geqslant q'(t)-vq(t)\\[0.2cm]
&=0,
\end{align*}
where we need $\gamma'(t)>0$.

On the other hand, when $u_1(x,t)\in[\theta,1-\theta]$, due to $\partial_tu_1(x,t)>0$, there is a constant $\overline{b}>0$ such that $\partial_tu_1(x,t)\geqslant \overline{b}$. Thus it follows from \eqref{eq:con} and \eqref{eq:const} that
\begin{align*}
\partial_t\overline{u}_1^{+}-\partial_{xx}\overline{u}_1^{+}
-f(\overline{u}_1^{+})&=\gamma'(t)
\partial_tu_1+\partial_tu_1-\partial_{xx}u_1+q'(t)-f(u_1+q(t))\\[0.2cm]
&\geqslant \overline{b}\gamma'(t)+vq(t)-wq(t)\\[0.2cm]
&=0,
\end{align*}
where we need the function $\gamma(t)$ satisfying
\begin{equation}\label{eq:ode12}\left\{\begin{array}{ll}
\gamma'(t)=\frac{w-v}{\overline{b}}\overline{q}_0e^{vt},\quad t>0,\\[0.2cm]
\gamma(0)=\overline{q}_0.
\end{array}
\right.
\end{equation}
Solving the equations \eqref{eq:ode1} and \eqref{eq:ode12} yields that $q(t)=\overline{q}_0e^{vt}$ and $\gamma(t)=\overline{q}_0(1+\gamma_0
-\gamma_0e^{vt})$, where $\gamma_0=\frac{v-w}{\overline{b}v}>0$, since $v<0$ and $w>0$. Also since $\partial_tu_1(x,t)>0$, then $u_1(x,0)\leqslant u_1(x,\overline{q}_0)+\overline{q}_0=\overline{u}_1^{+}(x,0)$. Hence,
$\overline{u}_1^{+}(x,t)$ is the supersolution of \eqref{eq:bi} with the initial condition $u_0(x)=u_1(x,0)$. Similarly, one can prove that $\underline{u}_1^{-}(x,t)$ is the subsolution of \eqref{eq:bi}  with the initial condition $u_0(x)=u_1(x,0)$.

Now we prove the local stability of the entire solution $u_1(x,t)$. For any given $\epsilon>0$, from Lemma \ref{thm:sest}, there exists a positive constant $\delta_1\leqslant\frac{\epsilon}{2L}$, such that for any $|t_3|\leqslant\delta_1$ and $t\in(1+t_3,+\infty)$,
\begin{equation}\label{eq:2est}
\|u_1(\cdot,t+t_3)-u_1(\cdot,t)\|=\|\partial_tu_1
(\cdot,t+t^*)\||t_3|
\leqslant\frac{\epsilon}{2},
\end{equation}
where $t^*\in(-t_3,t_3)$. Choose $\overline{q}_0=\delta\leqslant\min\{\frac{\delta_1}{1+\gamma_0},
\frac{\epsilon}{2}\}$, and the initial function $u_0(x)$ satisfies $\|u_0(x)-u_1(x,0)\|<\delta$. Since $\partial_tu_1(x,t)>0$, then $u_1(x,0)+\delta\leqslant u_1(x,\delta)+\delta$, $u_1(x,-\delta)-\delta\leqslant u_1(x,0)-\delta$. Therefore, for all $x\in\mathbb{R}$, we have
\begin{align*}
u_1(x,-\delta)-\delta\leqslant u_1(x,0)-\delta\leqslant u_0(x)\leqslant u_1(x,0)+\delta\leqslant u_1(x,\delta)+\delta.
\end{align*}
Thus $\overline{u}_1^{+}(x,t)$ and $\underline{u}_1^{-}(x,t)$ are also the supersolution and subsolution of \eqref{eq:bi} with the initial condition $u_0(x)$. Consequently,
\begin{align}\label{eq:inequ}&\max\{0,u_1(x,t-\delta(1+\gamma_0-\gamma_0e^{vt}))
-\delta e^{vt}\}\nonumber\\[0.2cm]\leqslant& u(x,t;u_0)
\leqslant\min\{1,u_1(x,t+\delta(1+\gamma_0-\gamma_0e^{vt}))+\delta e^{vt}\}.
\end{align}
By noting that $\delta(1+\gamma_0-\gamma_0e^{vt})\leqslant\delta_1$, it follows from \eqref{eq:2est} that
$$u_1(x,t)-\epsilon\leqslant u(x,t;u_0)\leqslant u_1(x,t)+\epsilon.$$
In a word, for any $\epsilon>0$, there exists a $\delta>0$, when $\|u_0(x)-u_1(x,0)\|<\delta$, then
$$\|u(x,t;u_0)-u_1(x,t)\|\leqslant\epsilon, \ \ \ \ \ \ t\geqslant0,$$
which means that the entire solution $u_1(x,t)$ is local asymptotic stable.

Finally, we will show that the entire solution $u_1(x,t)$ is local exponential asymptotic stable.  First of all, since $\underline{u}_1(x,t)$ given in Theorem \ref{thm:entire1} is also the subsolution of \eqref{eq:bi} for $t\geqslant0$, then it follows from the comparison theorem that for all $(x,t)\in\mathbb{R}\times\mathbb{R}$,
\begin{equation}\label{eq:ineq}
\max\{\phi(x+ct+x_6),\phi(-x+ct+x_6)\}=\underline{u}_1(x,t)\leqslant u_1(x,t)<1.
\end{equation}
Therefore, from \eqref{eq:esti} and \eqref{eq:ineq}, there is a $T_3>0$ with $cT_3+x_6>0$ such that for any $t\geqslant T_3$,  we have
\begin{align}
0\leqslant1-u_1(x,t)\leqslant1-\underline{u}_1(x,t)
\leqslant M_4e^{\mu_2(|x|+ct+x_6)}
\leqslant M_4e^{\mu_2ct}.\label{eq:exp1}
\end{align}

The next step is to prove
\begin{equation}\label{eq:sub}\underline{u}_1^{-}(x,t)=u_1(x,t-\delta(1+\gamma_0-\gamma_0e^{vt}))
-\delta e^{vt}
\end{equation}
holds for large $t$. Indeed, we remark that
\begin{equation}\label{eq:msub}
\partial_t(u_1(x,t-\overline{q}_0(1+\gamma_0-\gamma_0e^{vt}))
-\overline{q}_0e^{vt})=(1+\overline{q}_0\gamma_0ve^{vt})
\partial_tu_1-\overline{q}_0ve^{vt}.
\end{equation}
Since $\gamma_0=\frac{v-w}{\overline{b}v}$, then $$1+\overline{q}_0\gamma_0ve^{vt}=1-\overline{q}_0
\frac{w-v}{\overline{b}}e^{vt}.$$
Hence we can choose a constant $T_4\geqslant\max\{0,\frac{1}{v}
\ln\frac{\overline{b}}{(w-v)\theta}\}$ such that when $t\geqslant T_4$, by noting that $w>0$, $v<0$, $\overline{b}>0$ as well as $0\leqslant\overline{q}_0\leqslant\theta$, we have $$\overline{q}_0\frac{w-v}{\overline{b}}e^{vt}\leqslant\theta
\frac{w-v}{\overline{b}}e^{vt}\leqslant\theta
\frac{w-v}{\overline{b}}e^{vT_4}\leqslant1.$$
Therefore, when $t\geqslant T_4$, from $v<0$ and $\partial_tu_1(x,t)>0$ as well as \eqref{eq:msub}, we can obtain $$\partial_t(u_1(x,t-\overline{q}_0(1+\gamma_0-\gamma_0e^{vt}))
-\overline{q}_0e^{vt})>0.$$
Moreover, by noting $\lim\limits_{t\rightarrow+\infty}\underline{u}_1
(x,t-\overline{q}_0(1+\gamma_0-\gamma_0e^{vt}))=1$, thus,
there is a constant $T_5\geqslant\max\{T_3,\ T_4,\ \delta
(1+\gamma_0),\ \delta
(1+\gamma_0)-\frac{x_6}{c}\}$ such that for $t\geqslant T_5$, \eqref{eq:sub} holds.

Therefore, for $t\geqslant T_5\geqslant \delta(1+\gamma_0)$, it follows from \eqref{eq:inequ} and \eqref{eq:sub} that
\begin{align}
|1-u(x,t;u_0)|&\leqslant|1-u_1(x,t-\delta(1+\gamma_0-
\gamma_0e^{vt}))+\delta e^{vt}|\nonumber\\[0.2cm]
&\leqslant M_4e^{\mu_2(|x|+ct-c\delta(1+\gamma_0-\gamma_0e^{vt})+x_6)}
+\delta e^{vt}\nonumber\\[0.2cm]
&\leqslant M_4e^{\mu_2ct}+\delta e^{vt}.\label{eq:exp3}
\end{align}
Thus, when $\|u_0(x)-u_1(x,0)\|<\delta$,  for $t\geqslant T_5$, from \eqref{eq:exp1} and \eqref{eq:exp3}, we have
\begin{align*}
\|u(x,t;u_0)-u_1(x,t)\|
\leqslant&\|1-u_1(x,t)\|
+\|1-u(x,t;u_0)\|\\[0.2cm]
\leqslant& 2M_4e^{\mu_2ct}+
\delta e^{vt}.
\end{align*}

Next we discuss the stability of $u_1(x,t)$ under the case (C2). The proof of the local stability is similar, we only need to consider the local exponential stability. By directly calculating \eqref{eq:2ode}, we know that
$$p_3(t)=p_3(0)+ct-\frac{1}{\lambda_1}\ln\left\{1-\frac{M_8}{c}
e^{\lambda_1p_3(0)}
(1-e^{c\lambda_1t})\right\}.$$
Set $$g(t)=\frac{1}{\lambda_1}\ln\left\{1-\frac{M_8}{c}e^{\lambda_1p_3(0)}
(1-e^{c\lambda_1t})\right\}.$$
Then $$g'(t)=\frac{cM_8e^{\lambda_1p_3(0)}e^{\lambda_1t}}{c-M_8
e^{\lambda_1p_3(0)}
+M_8e^{\lambda_1p_3(0)}e^{c\lambda_1t}}>0,$$
since $p_3(0)<\frac{1}{\lambda_1}\ln\frac{c}{M_8}$.
Thus, $g(t)<g(\infty)=\frac{1}{\lambda_1}\ln\{1-\frac{M_8}{c}
e^{\lambda_1p_3(0)}\}:=x_7$.
Consequently, $p_3(t)>ct+x_8$, where $x_8:=p_3(0)-x_7$. As a result, $$\underline{u}_1(x,t)>\max\{\phi(x+ct+x_8),\phi(-x+ct+x_8)\}.$$
Thus in the case (C2), we also get an inequality similar to \eqref{eq:ineq}. The rest of the proof is similar.
\qed

\begin{remark}\label{remark:remark1}
Under the cases (C1) and (C2), by noting $\int^1_0f(u)du>0$,  the entire solution $u_1(x,t)$ found in Theorem \ref{thm:entire1} satisfies $$\lim\limits_{t\rightarrow+\infty}\|u_1(x,t)-1\|=0.$$ In fact, this conclusion coincides with the results in \cite{y03}. Thus the super-sub solution method is a valid way to simplify the proof of the existence of entire solutions of \eqref{eq:bi} in \cite{y03}.
\end{remark}

\begin{remark}\label{remark:remark2}
Under the case (C1), since the subsolution is $$\underline{u}_1(x,t)=\max\{\phi(x+ct+x_6),\phi(-x+ct+x_6)\},$$ then the existence, uniqueness and stability of entire solutions of \eqref{eq:bi} can be found in \cite{wlr09} as well. Here in order to prove the local asymptotic stability of entire solutions of \eqref{eq:bi}, we constructed the different supersolution and subsolution compared with \cite{wlr09}. Moreover we also simplify the way to prove the local asymptotic stability of entire solutions of \eqref{eq:bi} compared with \cite{y03}.
\end{remark}

\begin{remark}\label{remark:remark3}
When $f(u)=u(1-u)(u-\alpha)\ \alpha\in(0,1)$, one easily sees that $\int^1_0f(s)ds=\frac{1-2\alpha}{12}$ and $f'(0)=-\alpha$, $f'(1)=-1+\alpha$. Obviously,
$$f'(0)\lesseqqgtr f'(1)\ if\ and\ only\ if\ \int^1_0f(u)du\lesseqqgtr0\ if\ and\ only\ if\ \alpha\gtreqqless\frac{1}{2}.$$
Therefore if $\int^1_0f(u)du>0$, then $f'(0)>f'(1)$, namely, only the case (C1) will occur.
\end{remark}

For the case $\int^1_0f(u)du<0$, the authors in \cite{cg05} had found the entire solution $u_2(x,t)$ for \eqref{eq:bi}. Now we will discuss the long time behavior and the local exponential asymptotic stability of the entire solution in the following theorem.

\begin{theorem}\label{thm:entire3}
Assume that $\int^1_0f(u)du<0$ and $\phi$ is the solution to \eqref{eq:obi}, then \eqref{eq:bi} admits a unique entire solution $u_2(x,t)$ satisfying $\partial_tu_2(x,t)<0$, $u_2(x,t)=u_2(-x,t)$, $0<u_2(x,t)<1$,
and for $(x,t)\in\mathbb{R}\times(-\infty,-4B\phi(0)],$
$$u_2(x,t+h_1(t))<\phi(-x+ct)\phi(x+ct)<u_2(x,t-h_1(t)),$$
where $h_1(t)=4B\phi(ct)$ and $B>0$, 
$\lim\limits_{t\rightarrow+\infty}\|u_2(x,t)\|=0$,
\begin{equation}\label{eq:asy1}
\lim\limits_{t\rightarrow-\infty}\{\sup\limits_{x\geqslant0}
|u_2(x,t)-\phi(-x+ct)|
+\sup\limits_{x\leqslant0}|u_2(x,t)-\phi(x+ct)|\}=0.
\end{equation}
Moreover, the unique entire solution $u_2(x,t)$ is local exponential asymptotic stable.
\end{theorem}
\Proof Firstly, it follows from \cite{cg05} that
\begin{equation}\label{eq:asye}
\lim\limits_{t\rightarrow-\infty}\|u_2(x,t)
-\phi(x+ct)\phi(-x+ct)\|=0,
\end{equation}
which implies that
\begin{align*}
\lim\limits_{t\rightarrow-\infty}\sup\limits_{x\geqslant0}|u_2(x,t)-\phi(-x+ct)|
\leqslant&\lim\limits_{t\rightarrow-\infty}\sup\limits_{x\geqslant0}
|u_2(x,t)-\phi(-x+ct)\phi(x+ct)|\\[0.2cm]
+&\lim\limits_{t\rightarrow-\infty}\sup\limits_{x\geqslant0}|\phi(-x+ct)
||1-\phi(x+ct)|\\[0.2cm]
=&0.
\end{align*}
Similarly, $\lim\limits_{t\rightarrow-\infty}\sup\limits_{x\leqslant0}|u_2(x,t)-\phi(x+ct)|=0$. In a word, the entire solution $u_2(x,t)$ satisfies \eqref{eq:asy1}.

Secondly, since $\int^1_0f(u)du<0$, then the wave speed $c<0$, and at the same time, according to \cite{cg05}, the entire solution $u_2(x,t)$ satisfies for $t\leqslant-4B\phi(0)<0$, $$u_2(x,t+h_1(t))<\phi(-x+ct)\phi(x+ct)<u(x,t-h_1(t)).$$
Hence, for any $\tau\leqslant-4B\phi(0)<0$, we have
$$u_2(x,\tau+h_1(\tau))<\phi(x+c\tau),$$
and the comparison theorem yields that for all $t>\tau$, $u_2(x,t+h_1(\tau))<\phi(x+ct)$. Since $h_1(t)=4B\phi(ct)$, setting $\tau\rightarrow-\infty$ leads to $u_2(x,t)\leq\phi(x+ct)$. Similarly, $u_2(x,t)\leq\phi(-x+ct)$. Thus, for all $(x,t)\in\mathbb{R}
\times\mathbb{R}$,
\begin{equation}\label{eq:2est1}
u_2(x,t)\leq\min\{\phi(x+ct),\phi(-x+ct)\}.
\end{equation}
Specially, $u_2(x,0)\leq\min\{\phi(x),\phi(-x)\}$, which implies that $\lim\limits_{x\rightarrow\pm\infty}u_2(x,0)=0<\alpha$. Hence, due to Lemma \ref{lem:stacon}, $\lim\limits_{t\rightarrow+\infty}\|u_2(x,t)\|=0$.

Thirdly, similar to the proof in Theorem \ref{thm:entire1}, since $f$ satisfies the assumption (A), there are some constants $\alpha_1<\alpha<\alpha_2$ such that $f\neq0$ in $(0,\alpha_1]\cup[\alpha_2,1)$. Then from the monotonicity of $\phi$, there exists a positive constant $\hat{x}$ such that $\phi(-\hat{x})\leqslant\alpha_1$.
Set $$l_2(t)=-ct-\hat{x}, \ \ \ \ m_1(t)=ct+\hat{x},$$ then for any $t\leqslant0$, obviously $m_2(t)\geqslant0\geqslant l_2(t)$, and according to \eqref{eq:2est1}, $u_2(x,t)\leqslant\alpha_1$ for any $x\in(-\infty,l_2(t)]\cup[m_2(t),+\infty)$. On the other hand, according to \eqref{eq:asy1}, for any $\varepsilon>0$, there exists a $T_{13}<0$ such that for any $t\leqslant T_{13}$,
$$\sup\limits_{x\geqslant0}|u_2(x,t)-\phi(-x+ct)|
\leqslant\varepsilon,\quad\sup\limits_{x\leqslant0}|u_2(x,t)-\phi(x+ct)|
\leqslant\varepsilon,$$
which together with the monotonicity of $\phi$ implies that
$$u_2(x,t)\geqslant\phi(-m_2(t)+d_2+ct)-\varepsilon\geqslant\alpha_2$$
holds for $x\leqslant-m_2(t)+d_2=-ct+\hat{x}+d_2$ with some sufficiently large positive constant $d_2$. Similarly, $u_2(x,t)\geqslant\alpha_2$ for $x\geqslant l_2(t)+d_2$. Since there is a $T_{14}<0$ satisfying $l_2(t)+d_2\leqslant0\leqslant m_2(t)-d_2$ for any $t\leq T_{14}$, if we choose  $T_2\leqslant\min\{T_{13},T_{14}\}$, then $u_2(x,t)\geqslant\alpha_2$ for any $t\leqslant T_2$ and $x\in[l_2(t)+d_2,m_2(t)-d_2]$. Therefore, the entire solution $u_2(x,t)$ satisfies the condition $\mathbb{M}^-$, and from the proof in the papers \cite{cg05} and \cite{wlr09}, it is unique except for a space-time translation.

Fourthly, we consider the local asymptotic exponential stability of $u_2(x,t)$. First of all, for all $(x,t)\in\mathbb{R}\times\mathbb{R}$, similar to the proof of Theorem \ref{thm:entire1l}, it is easy to prove $\partial_tu_2(x,t)<0$. Thus there exists a negative constant $\tilde{b}$ such that when $u_2(x,t)\in[\theta,1-\theta]$, $\partial_tu_2(x,t)\leqslant\tilde{b}<0$. Similarly, it is not hard to prove the following two functions
\begin{align*}
\overline{u}_2^{+}(x,t)=\min\{1,u_2(x,t-\tilde{\gamma}(t))+q(t)\}\\[0.2cm] \underline{u}_2^{-}(x,t)=\max\{0,u_2(x,t+\tilde{\gamma}(t))-q(t)\}
\end{align*}
are the supersolution and subsolution of \eqref{eq:bi} with the initial condition $u_0(x)=u_2(x,0)$, respectively. At this time, the function $q(t)$ also satisfies \eqref{eq:ode1}, while the function $\tilde{\gamma}(t)$ satisfies
\begin{equation*}\left\{\begin{array}{ll}
\tilde{\gamma}'(t)=\frac{w-v}{-\tilde{b}}\overline{q}_0e^{vt},\quad t>0,\\[0.2cm]
\tilde{\gamma}(0)=\overline{q}_0,
\end{array}
\right.
\end{equation*}
where the parameters $\overline{q}_0$, $\theta$ are defined in the above, as well as $w$ and $v$ are defined in \eqref{eq:con} and \eqref{eq:const} separately. Then from the proof in Theorem \ref{thm:entire1l} or \cite{wlr09}, we know that $u_2(x,t)$ is local Lyapunov stable. In the end, we discuss the local exponential asymptotic stability of $u_2(x,t)$. With the similar proof, it is not hard to see that there is a constant $T_6$ such that for $t\geqslant T_6$, $\overline{u}_2^{+}(x,t)=u_2(x,t-\tilde{\gamma}(t))+q(t)$. By choosing $T_7=\max\{T_6,\delta(1+\tilde{\gamma}_0)\}$, where $\tilde{\gamma}_0=\frac{w-v}{\tilde{b}v}>0$, and noting $c<0$, then when $t\geqslant T_7$,
 with the help of \eqref{eq:esti} and \eqref{eq:2est1}, we finally arrive at
\begin{align*}
\|u(x,t;u_0)-u_2(x,t)\|
\leqslant&\|u_2(x,t)\|
+\|u(x,t;u_0)\|\\[0.2cm]
\leqslant& u_2(x,t)+u_2(x,t-\delta(1+\tilde{\gamma}_0-\tilde{\gamma}_0e^{vt}))
+\delta e^{vt}\\[0.2cm]
\leqslant& 2M_3e^{\lambda_1ct}+
\delta e^{vt},
\end{align*}
which implies the local exponential asymptotic stability of $u_2(x,t)$. Thus we have finished the proof.\qed


\section{The monostable equation}
In this section, under the assumption (A$'$), we will discuss the local exponential asymptotic stability of entire solutions of \eqref{eq:bi}.

From \cite{gm05} there is an entire solution $u_3(x,t)$ satisfying
\begin{align*}
&\max\{\phi_{c_1}(x+c_1t+y_3),\phi_{c_2}(-x+c_2t+y_4)\}\\[0.2cm]
\leqslant &u_3(x,t)
\leqslant\min\{1,\phi_{c_1}(x+p_4(t))
+\phi_{c_2}(-x+p_5(t))\},
\end{align*}
where $\phi_{c_k}\ (k=1,2)$ are the solutions to \eqref{eq:obi}, $c_1,c_2\in[c_{min},+\infty)$, $c_{min}=2\sqrt{f'(0)}$, $c_1\leqslant c_2$, and the functions $p_4(t)$, $p_5(t)$ are the solutions to
\begin{align*}\left\{\begin{array}{ll}
p'_4(t)=c_1+M_9e^{\tilde{\alpha}p_4},\\[0.2cm]
p'_5(t)=c_2+M_9e^{\tilde{\alpha}p_4},
\end{array}\right.
\end{align*}
where $p_5(0)\leqslant p_4(0)\leqslant0$, and $M_9$, $\tilde{\alpha}$ are positive constants.

Suppose that the function $\rho(t)$ is the solution to $\rho'=f(\rho)$ with $0<\rho(t)<1$. Define $\nu(t):=\nu_0e^{f'(0)t}\ (\nu_0>0)$ such that $0<\rho(t)-\nu(t)\leqslant M_{10}e^{f'(0)t}$ for $t\leqslant0$. From \cite{gm05}, there are three entire solutions $u_{ij}(x,t)$, $(i,j)=(1,0),(0,1),(1,1)$ satisfying
\begin{align*}
&\max\{\chi_i\phi_{c_1}(x+c_1t+y_6),\chi_j\phi_{c_2}(-x+c_2t+y_7),
\rho(t)\}\\[0.2cm]\leqslant& u_{ij}(x,t)
\leqslant\min\{1,\chi_i\phi_{c_1}(x+p_4(t))
+\chi_j\phi_{c_2}(-x+p_5(t))+\nu(t)\},
\end{align*}
where $\chi_i=i$.

We remark that whenever $\underline{u}_3(x,t)=\max\{\phi_{c_1}(x+c_1t+y_4),
\phi_{c_2}(-x+c_2t+y_5)\}$, or $\max\{\chi_i\phi_{c_1}(x+c_1t+y_6),\chi_j\phi_{c_2}(-x+c_2t+y_7),
\rho(t)\}$, then $\partial_t\underline{u}_3(x,t)\geqslant0$. Then similar to the proof in Theorem \ref{thm:entire1l}, it is not hard to see that $\partial_tu_{ij}(x,t)>0\ ((i,j)=(1,0),(0,1),(1,1))$ and $\partial_tu_3(x,t)>0$ for all $(x,t)\in\mathbb{R}\times\mathbb{R}$.
Here, the proofs of the local stability and the local exponential asymptotic stability of the entire solutions are similar to those in the bistable case. By taking the entire solution $u_3(x,t)$ for example, we only point out the different part. Since $f'(1)<0$ and $f'(0)>0$, there are positive constants $\theta'$ and $\theta''$ such that $f'(u)<0\ (u\in[1-2\theta', 1+\theta'])$ and $f'(u)>0\ (u\in[0, 2\theta''])$. For convenience, we still set $\theta=\min\{\theta',\theta''\}$. Because $\partial_tu_3(x,t)>0$, there exists some constant $T_8$ such that for all $(x,t)\in\mathbb{R}\times(T_8,+\infty)$, $u_3(x,t)>\theta$. Now, we consider the initial problem starting at time $\max\{0,T_8\}$, and then similar to the proof of Theorem \ref{thm:entire1l}, $u_3(x,t)$ is local exponential asymptotic stable. In a word, we obtain the following theorem.
\begin{theorem}\label{thm:entire5}
Suppose that $f$ satisfies (A$'$) and $f'(u)\leqslant f'(0)\ \ (u\in[0,1])$. Furthermore, let $\phi_{c_k}$ be the solutions to \eqref{eq:obi} with $c_k\in[c_{min},+\infty)$, $k=1,2$. Then the entire solutions $u_3$ and $u_{ij}(x,t)$, $(i,j)=(1,0),(0,1),(1,1)$ are local exponential asymptotic stable.
\end{theorem}

\section*{References}
\bibliographystyle{elsarticle-num}

\begin{thebibliography}{aa}
\bibitem{al91} S. Ahmad, A. C. Lazer, An elementary approach to traveling front solutions to a system of $N$ competition-diffusion equations, Nonlinear Anal. 16 (1991) 892-901.

\bibitem{alt08} S. Ahmad, A. C. Lazer, A. Tineo, Traveling waves for a system of equations, Nonlinear Anal. 68 (2008) 3909-3912.

\bibitem{aw75} D. G. Aronson, H. F. Weinberger, Nonlinear diffusion in population genetics, combustion, and nerve pulse propagation. Partial Differential Equations and Related Topics, Lecture Notes in Math., vol. 446, Springer, Berlin, 1975, pp.5-49.

\bibitem{aw78} D. G. Aronson, H. F. Weinberger, Multidimensional nonlinear diffusion arising in population genetics, Adv. in Math. 30 (1978) 33-76.

\bibitem{cg05} X. Chen, J. S. Guo, Existence and uniqueness of entire solutions for a reaction-diffusion equation, J. Differential Equations 212 (2005) 62-84.

\bibitem{c71} H. Cohen, Nonlinear diffusion problems, Studies in Mathematics 7, Studies in Applied Mathematics, ed. A. Taubo Math. Assoc. of America and Prentice Hall, (Englewood Cliffs, N. J.) 1971, 27-63.

\bibitem{cg84} C. Conley, R. Gardner,  Application of the generalized Morse index to travelling wave solutions of a competitive reaction-diffusion model, Indiana Univ. Math. J. 33 (1984) 319-343.

\bibitem{ct12} E. C. M. Crooks, J. C. Tsai, Front-like entire solutions for equations with convection, J. Differential Equations 253 (2012) 1206-1249.

\bibitem{f79l} P. C. Fife, Mathematical aspects of reacting and diffusing systems, Lecture Notes in Biomathematics, vol. 28, Springer-Verlag, Berlin, 1979.

\bibitem{f79} P. C. Fife, Long time behavior of solutions of bistable nonlinear diffusion equations, Arch. Ration. Mech. Anal. 70 (1979) 31-46.

\bibitem{fm77} P. C. Fife, J. B. McLeod, The approach of solutions of nonlinear diffusion equations to travelling front solutions, Arch. Ration. Mech. Anal. 65 (1977) 335-361.

\bibitem{fm81} P. C. Fife, J. B. McLeod, A phase plane discussion of convergence to travelling fronts for nonlinear diffusion, Arch. Ration. Mech. Anal. 75 (1981) 281-314.

\bibitem{f37} R. A. Fisher, The advance of advantageous genes, Ann. of Eugenics 7 (1937) 355-369.

\bibitem{fmn04} Y. Fukao, Y. Morita, H. Ninomiya, Some entire solutions of the Allen-Cahn equation, Taiwanese J. Math. 8 (2004) 15-32.

\bibitem{gr07} T. Gallay, E. Risler, A variational proof of global stability for bistable travelling waves, Differential Integral Equations 20 (2007) 901-926.

\bibitem{g82} R. Gardner, Existence and stability of travelling wave solutions of competition models: a degree theoretic approach, J. Differential Equations 44 (1982) 343-364.

\bibitem{gl13} J. S. Guo, Y. C. Lin, Entire solutions for a discrete diffusive equation with bistable convolution type nonlinearity, Osaka J. Math. 50 (2013) 607-629.

\bibitem{gm05} J. S. Guo, Y. Morita, Entire solutions of reaction-diffusion equations and an application to discrete diffusive equations, Disc. Cont. Dyn. Syst. 12 (2005) 193-212.

\bibitem{hn99} F. Hamel, N. Nadirashvili, Entire solutions of the KPP Equation, Comm. Pure Appl. Math. 52 (1999) 1255-1276.

\bibitem{hn01} F. Hamel, N. Nadirashvili, Travelling fronts and entire solutions of the Fisher-KPP Equation in R$^N$, Arch. Rational Mech. Anal. 157 (2001) 91-163.

\bibitem{k96} Y. Kan-on, Existence of standing waves for competition-diffusion equations, Japan J. Indust. Appl. Math., 13 (1996) 117-133.

\bibitem{kpp37} A. Kolmogoroff, I. Petrovsky, N. Piscounoff, \'{E}tude de \'{I}equation de la diffusion avec croissance
de la quantit\'{e} de mati\`{e}re et son application \`{a} unprobleme biologique, Bull. Univ. Moskou, Ser. Internat., Sec. A, 1 (1937) 6, 1-25.

\bibitem{lx92} C. D. Levermore, J. X. Xin, Multidimensional stability of traveling waves in a bistable reaction-diffusion equation. II, Comm. Partial Differential Equations 17 (1992) 1901-1924.

\bibitem{mn06} Y. Morita, H. Ninomiya, Entire solutions with merging fronts to reaction-diffusion equations, J. Dynam. Differential Equations 18 (2006) 841-861.

\bibitem{mt09} Y. Morita, K. Tachibana, An entire solution to the Lotka-Volterra competition-diffusion equations, SIAM J. Math. Anal. 40 (2009) 2217-2240.

\bibitem{nya65} J. Nagumo, S. Yoshizawa, S. Arimoto, Bistable transmission lines, I.E.E.E. Transactions on Circuit Theory 12 (1965) 400-412.

\bibitem{om99d} T. Ogiwara, H. Matano, Monotonicity and convergence results in order preserving systems in the presence of symmetry, Disc. Cont. Dyn. Syst. 5 (1999) 1-34.

\bibitem{om99p} T. Ogiwara, H. Matano, Stability analysis in order-preserving systems in the presence of symmetry, Proc. Roy. Soc. Edinburgh Sect. A 129 (1999) 395-438.

\bibitem{r08} E. Risler, Global convergence toward traveling fronts in nonlinear parabolic systems with a gradient structure, Ann. Inst. H. Poincar\'{e} Anal. Non Lin\'{e}aire 25 (2008) 381-424.

\bibitem{r04} V. Roussier, Stability of radially symmetric travelling waves in reaction-diffusion equations, Ann. Inst. H. Poincar\'{e} Anal. Non Lin\'{e}aire 21 (2004) 341-379.

\bibitem{wl15} Y. Wang, X. Li, Some entire solutions to the competitive reaction diffusion system, J. Math. Anal. Appl. 430 (2015) 993-1008.

\bibitem{wlr09} Z. C. Wang, W. T. Li, S. Ruan, Entire solutions in bistable reaction-diffusion equations with nonlocal delayed nonlinearity, Trans. Amer. Math. Soc. 361 (2009) 2047-2084.

\bibitem{x92} J. X. Xin, Multidimensional stability of traveling waves in a bistable reaction-diffusion equation. I. Comm. Partial Differential Equations 17 (1992) 1889-1899.

\bibitem{y03} H. Yagisita, Backward global solutions characterizing annihilation dynamics of travelling fronts, Publ. Res. Inst. Math. Sci 39 (2003) 117-164.

\end{thebibliography}

\end{document}